\documentstyle[12pt]{article}
\newcommand{\const}{\mathop{\rm const}\limits}

\newcommand{\supp}{\mathop{\rm supp}\limits}

\newcommand{\card}{\mathop{\rm card}\limits}

\newcommand{\Var}{\mathop{\rm Var}\limits}

\newcommand{\Law}{\mathop{\rm Law}\limits}

\newcommand{\Sub}{\mathop{\rm Sub}\limits}

\newcommand{\SSub}{\mathop{\rm SSub}\limits}

\newcommand{\grad}{\mathop{\rm grad}\limits}

\begin{document}

\begin{center}

{\bf  BANACH SPACES CHARACTERIZATION OF RANDOM VECTORS \\

\vspace{3mm}

 WITH EXPONENTIAL DECREASING TAILS OF DISTRIBUTION. }\\

\vspace{5mm}

 {\sc Ostrovsky E., Sirota L.}\\

\vspace{3mm}

 \ Department of Mathematics and Statistics, Bar-Ilan University, \\
59200, Ramat Gan, Israel.\\
e-mail: \ eugostrovsky@list.ru \\

\vspace{3mm}

 \ Department of Mathematics and Statistics, Bar-Ilan University,\\
59200, Ramat Gan, Israel.\\
e-mail: \ sirota3@bezeqint.net \\

\vspace{3mm}

 {\sc Abstract.}

\end{center}

 \  We present in this paper the Banach space representation for the set of random finite-dimensional {\it vectors}
with exponential decreasing tails of distributions. \par
 \  We show that there are at last three types of these multidimensional Banach spaces, i.e. which can completely describe
the random vectors with  exponential decreasing tails of distributions: exponential Orlicz spaces, Young spaces and Grand
Lebesgue spaces. \par
 \ We  discuss in the last section the possible applications of obtained results. \par
 \vspace{3mm}

 {\it Key words and phrases:}  Random variable (r.v.) and random vector (r.v.), centered (mean zero) r.v., moment generating function,
rearrangement invariant Banach space of vector random variables (vectors), ordinary and exponential moments, subgaussian and strictly
subgaussian random vectors and variables, Grand Lebesgue Space (GLS) and norm, multivariate tail of distribution, Young-Orlicz function,
Tchebychev's-Markov's and Chernov's inequality, theorem of Fenchel-Moraux, Young inequality, upper and lower
estimates, non-asymptotical exponential estimations, Kramer's condition, saddle-point method.

\vspace{4mm}


\vspace{5mm}

\section{ Introduction.  Notations. Statement of problem. Previous results.} \par

\vspace{4mm}

 \hspace{6mm} Let $ (\Omega,F,{\bf P} ) $ be a probability space, $ \Omega = \{\omega\}. $ \par

 \vspace{3mm}

 \ Denote by $  \epsilon  = \vec{\epsilon} = \{ \epsilon(1),   \epsilon(2), \ldots,   \epsilon(d) \} $ the non - random
$  d \ - $ dimensional  numerical vector, $  d = 1,2,3,\ldots, $ whose components take the values $  \pm 1 $ only. \par
 \ Set $  \vec{1} = (1,1,\ldots,1) \in R^d_+ $   and as usually $ \vec{0} = (0,0,\ldots,0) \in R^d  $ \par

 \ Denote by $  \Theta = \Theta(d) = \{ \ \vec{\epsilon} \ \} $ {\it collection} of all such a vectors.  Note that
 $  \card \Theta = 2^d $  and $ \vec{1} \in \Theta.$  \par

 \ Another denotations. The inequality $ a = \{a(i) \} =\vec{a} > \vec{b} = b = \{ b(i) \}, \ i = 1,2,\ldots,d $
for two $ d \ - $ dimensional numerical vectors imply the coordinate-wise comparison:
$  \forall i = 1,2,\ldots,d  \ \Rightarrow a(i) > b(i). $\par
 \ For $ \vec{\epsilon} \in \Theta(d) $  and vector $ \vec{x} $ we introduce the coordinatewise product
as a $ d  \ - $ dimensional vector of the form

$$
\vec{\epsilon} \otimes \vec{x} \stackrel{def}{=} \{ \epsilon(1) \ x(1), \ \epsilon(2) \ x(2), \ \ldots,  \epsilon(d) \ x(d)  \}.
$$

\vspace{3mm}

 \ {\bf Definition 1.1.}\par

\vspace{3mm}

 \ Let $ \xi = \vec{\xi} = (\xi(1), \xi(2), \ldots, \xi(d) )  $ be a {\it centered} (mean zero) random vector such that each its
component $  \xi(j) $  satisfies the famous Kramer's condition. The {\it natural function} $ \phi_{\xi}= \phi_{\xi}(\lambda), \
 \lambda = \vec{\lambda} = (\lambda(1), \lambda(2), \ldots, \lambda(d))  \in R^d  $
for the random vector $  \xi  $ is defined as follows:

$$
e^{\phi_{\xi}(\lambda)} \stackrel{def}{=} \max_{\vec{\epsilon}} \
{\bf E} \ e^{ \sum_{j=1}^d  \epsilon(j) \lambda(j) \xi(j)}=
$$

$$
\max_{\vec{\epsilon}} \
{\bf E} \exp \{\epsilon(1) \lambda(1) \xi(1) + \epsilon(2) \lambda(2)\xi(2)+ \ldots +   \epsilon(d) \lambda(d)\xi(d) \}, \eqno(1.1)
$$
where $  "\max" $ is calculated over all the combinations of signs $ \epsilon(j) =  \pm 1, $ on the other words,
$ \vec{\epsilon} \in \Theta. $ \par

 \ This concept is slightly different from the classical notion of the {\it Moment Generating Function, } see \cite{Feller1}. Indeed,
the last variable is calculated alike in (1.1), but without the operation $ \ \max_{\vec{\epsilon}} \ $  etc. \par

 \ The complete description of such a natural $  \phi_{\xi}(\cdot) $ multivariate functions, based on the Bernstein's theorem,
 represented in the preprint \cite{Ostrovsky5}. \par

\vspace{3mm}

 \ {\bf Definition 1.2.}\par

\vspace{3mm}

 \ The {\it tail function } for the random vector $ \vec{\xi} \ $
   $ U(\vec{\xi}, \vec{x}), \ \vec{x} = (x(1), x(2), \ldots, x(d)),   $ where all the coordinates $ x(j)  $
  of the deterministic vector $ \vec{x} $ are non-negative, is defined as follows.

$$
U(\vec{\xi}, \vec{x}) \stackrel{def}{=} \max_{\vec{\epsilon}}
{\bf P} \left( \cap_{j=1}^d  \{  \epsilon(j) \xi(j) > x(j) \}  \right) =
$$

$$
 \max_{\vec{\epsilon}}
{\bf P}(\epsilon(1) \xi(1) > x(1), \ \epsilon(2) \xi(2) > x(2), \ \ldots, \ \epsilon(d) \xi(d) > x(d) ), \eqno(1.2)
$$
 where as before $  "\max" $ is calculated over all the combinations of signs $ \epsilon(j) =  \pm 1. $\par

 \ We illustrate this notion at first in the case  $  d = 1:  $

$$
U(\xi,x) = \max( {\bf P}(\xi > x), \ {\bf P}(\xi < - x) ), \ x > 0.
$$

 \ Assume now $  d = 2.  $ Let $ \vec{\xi} = (\xi(1), \ \xi(2))  $ be a two-dimensional random
vector and let $ x, y  $  be non-negative numbers. Then

$$
U( (\xi(1), \xi(2)), \ (x,y) ) =
$$

$$
\max [ {\bf P} (\xi(1) > x, \ \xi(2) > y), \ {\bf P} (\xi(1) > x, \ \xi(2) < - y),
$$

$$
{\bf P} (\xi(1) < - x, \ \xi(2) > y), \ {\bf P} (\xi(1) < - x, \ \xi(2) < - y) ].
$$

\vspace{3mm}

 \ {\bf Definition 1.3.}\par

\vspace{3mm}

 \ Let $  h = h(x), \ x \in R^d  $ be some non-negative real valued function, which is finite at last on some non-empty
neighborhood of origin.  We denote as ordinary

$$
\supp h = \{x, \ h(x) < \infty   \}.
$$

 \ The Young-Fenchel, or Legendre transform $  h^*(y), \ y \in R^d  $ is defined likewise the one-dimensional case

$$
h^*(y) \stackrel{def}{=} \sup_{x \in \supp h} ( (x,y) - h(x)), \eqno(1.3)
$$

 \ Herewith $ (x,y)  $ denotes the inner (scalar) product of the vectors $ x,y: $

 $$
  \ (x,y) = \sum_{j=1}^d x(j)y(j); \hspace{5mm} \ |x|:= \sqrt{(x,x)}.
 $$

 \ Obviously, if the set $ \supp h $ is central symmetric, then the function $ h^*(y) $ is even.\par

 \ We recall here the famous theorem of Fenchel-Moraux: if the function $  h = h(x)  $ is convex  and continuous,
then $  h^{**}(y) = h(y). $  \par

\ Note in addition that if the function $  h = h(x), \ x \in R^d $ is radial (spherical) function:  $ h(x) = z(|x|), $
then its Young-Fenchel transform is also radial, as well. \par

\vspace{3mm}

 \ {\bf Definition 1.4.}\par

\vspace{3mm}

 \ Recall, see \cite{Rao1}, \cite{Rao2}, \cite{Krasnosel'skii1}, that the function $  x \to g(x), \ x \in R^d, \ g(x) \in R^1_+  $
 is named multivariate Young, or Young-Orlicz function, if it is even, convex, non-negative, twice continuous differentiable,
 finite on the whole space $  R^d,$ and such that

$$
g(x) = 0 \ \Leftrightarrow x = 0,
$$
and

$$
\det \frac{\partial^2 g(\vec{0})}{\partial x^2} > 0. \eqno(1.4)
$$

 We explain in detail:

$$
 \frac{\partial g}{\partial x} = \left\{  \frac{ \partial g}{ \partial x(j)} \right\} = \grad g,  \hspace{5mm}
\frac{\partial^2 g}{\partial x^2} = \left \{ \frac{\partial^2 g}{\partial x(k) \ \partial x(j)}  \right\}, \
j,k = 1,2,\ldots,d. \eqno(1.5)
$$

 \ We assume in addition finally

$$
\lim_{|x| \to \infty} \ \min_{ = 1,2, \ldots d} \  | \ \partial g(x)/\partial x(j) \ | = \infty. \eqno(1.6)
$$

 \ We will denote the set of all such a functions by $  Y = Y(R^d) $ and denote also by $ D  $ introduced before matrix

$$
D = D_g := \frac{1}{2} \left\{ \frac{\partial^2 g(0)}{\partial x(k) \ \partial x(l)} \right\}.
$$

 \ Evidently, the matrix $ D = D_g $ is non - negative definite,  write $  D = D_g \ge \ge 0.  $ \par

 \ Here and in what follows the relation $  A \ge \ge B  $  or equally $ B \le \le A $ between two square matrixes $  A  $  and $  B  $
 of the same size $  n \times n $ implies as ordinary

$$
\forall \lambda \in R^d \ \Rightarrow  (A \lambda, \lambda) \ge ( B \lambda, \lambda).
$$

\vspace{3mm}

 \ Let us now recall the definition and some properties of the
 Banach spaces consisting on the random vectors with exponential decreasing tails of distribution. The detail explanations
 and proofs reader may be find in the article \cite{Ostrovsky5}.  \par

\vspace{3mm}

 \ {\bf  Definition 1.5. } \par

\vspace{3mm}

 \ Let the Young function  $  \phi(\cdot)  $ be from the set $  Y = Y(R^d): \supp \phi = R^d. $\par

 \ We will say by definition likewise the one-dimensional case, see \cite{Kozachenko1}, \cite{Ostrovsky1}, chapter 1,
 that the {\it centered} random vector (r.v)
 $ \xi = \xi(\omega) = \vec{\xi} = \{\xi(1), \xi(2), \ldots, \xi(d) \} $ with values in the space $  R^d $ belongs to the
 space $ B(\phi), $ or equally Young space, write $ \vec{\xi}\in  B(\phi),  $  if there exists certain non-negative constant
 $ \tau \ge 0 $ such that

$$
\forall \lambda \in R^d \ \Rightarrow
\max_{\vec{\epsilon} \in \Theta } {\bf E} \exp \left( \sum_{j=1}^d \epsilon(j) \lambda(j) \xi(j) \right) \le
\exp[ \phi(\lambda \cdot \tau) ]. \eqno(1.7)
$$

 \ The minimal value $ \tau $ satisfying (1.7) for all the values $  \lambda \in R^d, $
is named by definition as a $ B(\phi) \ $ norm of the vector $ \xi, $ write

$$
||\xi||B(\phi)  \stackrel{def}{=}
$$

$$
 \inf \left\{ \tau, \ \tau > 0: \ \forall \lambda:  \ \lambda \in R^d \ \Rightarrow
 \max_{\vec{\epsilon}  \in \Theta }{\bf E}\exp \left( \sum_{j=1}^d \epsilon(j) \lambda(j) \xi(j) \right) \le
 \exp(\phi(\lambda \cdot \tau)) \right\}. \eqno(1.8)
$$

 \ Emerging in the relations \ (1.7-(1.8) \ the Young-Orlicz function $ \phi = \phi(\lambda), \ \in R^d $
is said to be {\it generating function } for the correspondent space $ B(\phi). $ \par

\vspace{3mm}

 \ For example, the  generating function  $ \phi(\lambda) = \phi_{\xi}(\lambda) $ for these spaces may be picked by the
 following  so-called {\it natural} way:

$$
\exp[ \phi_{\xi}(\lambda ) ] \stackrel{def}{=}
\max_{\vec{\epsilon}} {\bf E} \exp \left( \sum_{j=1}^d \epsilon(j) \lambda(j) \xi(j) \right),
 \eqno(1.9)
$$
if of course the random vector  $  \xi $ is centered and has an exponential tail of distribution. This imply  that
the natural function  $ \phi_{\xi}(\lambda ) $ is finite at last on some non-trivial central symmetry neighborhood of origin,
or equivalently  the mean zero random vector  $  \xi $ satisfies the multivariate Kramer's condition. \par

 \ Obviously, for the natural function $ \phi_{\xi}(\lambda )  $

$$
||\xi||B(\phi_{\xi}) = 1.
$$

 \ It is easily to see that this  choice of the generating function $ \phi_{\xi} $ is optimal, but in the practical using
often this function can not  be calculated in explicit view,  but there is a possibility to estimate its. \par

 \vspace{3mm}

 \ Note that the expression for the norm $ ||\xi||B(\phi)  $ dependent aside from the function $  \phi  $
 only on the distribution $  \Law(\xi). $ Thus, this norm and correspondent space  $  B(\phi)  $ are rearrangement invariant (symmetrical)
in the terminology of the classical book  \cite{Bennet1}, see chapters 1,2. \par

\vspace{3mm}

 \ It is proved in particular in the report \cite{Ostrovsky5} that: \\

{\sc  the space $ B(\phi) $ with respect to the norm $ || \cdot ||B(\phi) $ and
ordinary algebraic operations is a rearrangement invariant vector Banach space.} \par

 \vspace{3mm}

 \ Moreover, the proposition $ \vec{\xi} \in B(\phi)  $ may be completely adequate characterized through
the tail behavior $ U(\vec{\xi}, \vec{x}) $ as $ |x| \to \infty $ for the random vector $ \xi.  $

\vspace{3mm}

 {\bf We intend to significantly weaken  the conditions in \cite{Ostrovsky5} for tail estimates as well as for
moment estimates for random variables from these spaces. }

\vspace{3mm}

 \ {\bf Remark 1.1.} In the article of  Buldygin V.V. and Kozachenko Yu. V. \cite{Buldygin1} was
 considered a particular case when $ \phi(\lambda) = \phi^{(B)}(\lambda)  = 0.5(B\lambda,\lambda),  $ where
$  B  $  is non-degenerate positive definite symmetrical matrix, as a direct generalization of the one-dimensional one notion,
belonging to J.P.Kahane \cite{Kahane1}.  \par
 \ The correspondent random vector $ \vec{\xi}  $ was named in \cite{Buldygin1} as a subgaussian r.v. relative the
matrix $  B: $ \par

$$
\forall \lambda \in R^d \ \Rightarrow {\bf E} e^{(\lambda, \xi)} \le e^{0.5 (B \lambda, \lambda) \ ||\xi||^2 }.
$$

 \ We will write in this case  $  \xi \in \Sub(B)  $ or more precisely $ \Law \xi \in \Sub(B). $ \par

\vspace{3mm}

 \ {\bf Remark 1.2.}  Suppose the r.v. $  \xi  $ belongs to the space $  B(\phi);  $  then evidently
 $  {\bf E} \vec{\xi} = 0. $ \par
  \ Suppose in addition that it has there the unit norm, then

 $$
 \Var ( \vec{\xi}) \le \le   D_\phi; \hspace{4mm} \Leftrightarrow D_\phi \ge \ge \Var ( \vec{\xi}).
 $$

\vspace{3mm}

 \ It is interest to note that there exist many (mean zero) random vectors $  \eta = \vec{\eta} $ for which

$$
{\bf E} e^{(\lambda,\eta)} \le e^{0.5 \ (\Var(\eta) \lambda, \ \lambda) }, \ \lambda \in R^d,
$$
see e.g.  \cite{Buldygin2}, \cite{Buldygin3}, chapters 1,2; \cite{Ostrovsky1}, p.53. V.V.Buldygin and
Yu.V.Kozatchenko in \cite{Buldygin2} named these vectors {\it  strictly subgaussian;}  notation
$ \xi \in \SSub $ or equally $ \Law (\xi) \in \SSub. $ \par

 \ V.V.Buldygin and Yu.V.Kozatchenko find also some applications of these notions. \par

\vspace{4mm}

\section{ Main result: exponential tails behavior. } \par

\vspace{4mm}

 \hspace{5mm} We study in this section the connections between exponential decreasing tails behavior for random vectors and its
 norm in the introduced above spaces. \par

\vspace{3mm}

{\bf A. Upper estimate.} \par

\vspace{3mm}

  \ Statement of problem: let the $  B(\phi) \ -  $ norm of the non-zero random vector $ \xi \in R^d $ be a given; we can
suppose without loss of  generality $  ||\xi||B(\phi) = 1. $ We want to get the sharp  (if possible) upper bound for the tail
function $ U(\xi,x). $ \par
 \ Conversely: let for the centered random vector $ \xi = \vec{\xi} $  the its tail function $ U(\xi,x) $  (or its upper estimate)
be given; it is required to estimate its norm $  ||\xi||B(\phi)  $ in this space. \par

 \ In particular, it is required  to establish the one-to-one relation between behavior of the tail function $ U(\xi,x) $   as
$  \min_j x(j) \to \infty $ and belonging of the centered random vector $ \xi  $ to the space $  B(\phi). $ \par

 \ Note at first that when we derive the upper bounds for tail function $ U(\xi,x) $ through the $  ||\xi||B(\phi), $
 we do not need impose strong condition on the function $ \phi = \phi(\lambda). $
 \ Namely, let $  \phi = \phi(\lambda), \ \lambda \in V \subset R^d $ be arbitrary non-negative real valued function,
which is finite at last on some non-empty neighborhood $ V $ of origin. Suppose for given centered $  d  \ - $ dimensional
random vector $ \xi = \vec{\xi}  $

$$
{\bf E} e^{(\lambda,\xi)} \le e^{ \phi( \lambda) }, \ \lambda \in V.  \eqno(2.1)
$$

 \ On the other words, $ ||\xi||B(\phi) \le 1, $ despite the function $ \phi(\cdot) $ is not supposed in general case as Young-Orlicz
 function. Indeed,

$$
{\bf E} \exp( \lambda, \xi ) \le \exp(\phi(\lambda)), \ \lambda \in V. \eqno(2.1a)
$$

 \ One can to extend  formally this function on the whole space $  R^n:  $

$$
\phi(\lambda) := \infty, \ \lambda \not \in V.
$$

\vspace{4mm}

 \ {\bf Proposition 2.1.} Assume the (centered) random vector $ \xi $ satisfies the condition (2.1), or equally (2.2).
We state that for all the non-negative vector $ x = \vec{x}  $ there holds

 $$
 U(\vec{\xi}, \vec{x}) \le \exp \left( - \phi^*(\vec{x}) \right) \ - \eqno(2.2)
 $$
the multidimensional generalization of Chernov's  inequality. \par

\vspace{3mm}

 \ {\bf Proof.}  Let for definiteness $  \ \vec{x} > 0; \ $ the case when $ \exists k \ \Rightarrow x(k) \le 0  $ may be
considered analogously. \par
 \ We use the ordinary Tchebychev-Markov  inequality

$$
U(\vec{\xi}, \vec{x}) \le  \frac{ e^{ \phi( \lambda) }}{ e^{( \lambda,x )}  } =
e^{ -( \lambda,x ) - \phi( \lambda)}, \ \lambda > 0. \eqno(2.3)
$$

 \ Since the last inequality  (2.3) is true for arbitrary non-negative vector $  \lambda \in V,  $

$$
U(\vec{\xi}, \vec{x}) \le \inf_{\lambda >>0} e^{( \lambda,x ) - \phi( \lambda)} \le
$$

$$
\exp \left\{- \sup_{\lambda \in R^d} \exp [ \ ( \lambda,x ) - \phi( \lambda) \ ] \right\}  =
\exp \left( - \phi^*(\vec{x}) \right). \eqno(2.4)
$$

\vspace{3mm}

{\bf B. Lower estimate.} \par

\vspace{3mm}

 \ This case is more complicated and it requires more restrictions.
 Let us introduce  the following conditions on the Young-Orlicz functions.\par

\vspace{3mm}

 \ {\bf Definition 2.1. } We will say  that the {\it Young-Orlicz function} $  \phi = \phi(x), \ \in  R^d $
 satisfies a condition $ K(\gamma), $ write:  $  \phi(\cdot) \in K(\gamma), $
 iff there exists a positive number  $ \gamma \in (0,1)  $ such that the following integral is finite:

$$
I_{\gamma}(\phi) \stackrel{def}{=} \int_{R^d_+} \exp \left\{ \phi^*(\gamma \cdot x) -   \phi^*( x)   \right\} \
dx < \infty. \eqno(2.5)
$$

\vspace{3mm}

 \ {\bf Theorem 2.1.} Suppose that the Young-Orlicz  function $  \phi(\cdot)  $ satisfies  the condition (2.5):
 $ \phi(\cdot) \in K(\gamma).  $  Suppose also the mean zero random vector $  \xi = \vec{\xi} $ satisfies the
 condition (2.2) for all the non-negative deterministic vector $  \vec{x}: \ \forall j = 1,2,\ldots,d \ \Rightarrow x(j) > 0  $

 $$
 U(\vec{\xi}, \vec{x}) \le \exp \left( - \phi^*(\vec{x}) \right).    \eqno(2.6)
 $$
 \ We propose that r.v. $ \ \vec{\xi} \ $ belongs to the space $  B(\phi): \ \exists  \ C = C(\phi) \in (0,\infty), $

$$
{\bf E} e^{(\lambda,\xi)} \le e^{ \phi( C \cdot \lambda) }, \ \lambda \in R^d. \eqno(2.7)
$$

 \vspace{3mm}

 \ {\bf Proof}  is'nt likewise to the one-dimensional case, see  \cite{Kozachenko1}, \cite{Ostrovsky1}, p. 19-40.\par

  \ Note first of all that the estimate (2.7) is obviously satisfied   for the values $  \lambda = \vec{\lambda} $ with
appropriate positive constant $  C = C(\phi), $
for the  values $ \vec{\lambda}  $ belonging to a Euclidean unit ball of the space $   R^d: \ |\lambda| \le 1, $ since
the r.v. $  \xi $ is centered and has a very light (exponential decreasing) tail of distribution. It remains to consider
further only the case when $  |\lambda| \ge 1.$ \par
 \ Let for definiteness $  \vec{\lambda} >> \vec{1}. $  \ We have using integration  by parts

$$
 {\bf E} e^{(\lambda,\xi)}  \le 1 + \prod_{k=1}^d |\lambda_k| \cdot \int_{R^d_+} e^{(\lambda,x) - \phi^*(x) } \ dx \stackrel{def}{=} 1 +
  \prod_{k=1}^d |\lambda_k| \cdot  J_{\phi}(\lambda).
$$

 \ It is sufficient to investigate the main part of the last integral, namely

$$
J_{\phi}(\lambda) := \int_{R_+^d} e^{(\lambda,x) - \phi^*(x) } \ dx. \eqno(2.8)
$$

 \ The saddle-point method tell us that as $ |\lambda| \to \infty  $

$$
\log J_{\phi}(\lambda) \sim \sup_{x \in R^d} \left[ (\lambda,x) - \phi^*(x) \right] =
\phi^{**}(\lambda) = \phi(\lambda),
$$
by virtue of theorem Fenchel-Moraux. Therefore

$$
J_{\phi}(\lambda) \le \exp ( \phi(C \cdot \lambda) ), \ |\lambda| \ge 1.
$$

 \ To be more precisely, we attract the consequence from the direct definition of Young-Fenchel (Legendre) transform

$$
(\lambda,x) \le \phi^*(x) + \phi(\lambda),
$$
 the so - called  Young inequality.\par

 \ Let $  \gamma = \const \in (0,1).  $   We can write

$$
(\lambda,x) \le \phi^*(\gamma x) + \phi(\lambda /\gamma).
$$

  We have  after substituting into (2.8)

$$
J_{\phi}(\lambda) \le \exp \left(\phi(\lambda\gamma)  \right) \cdot \int_{R^d_+} \exp \left( \phi^*(\gamma x) - \phi(x)  \right) \ dx =
$$

$$
 I_{\gamma}(\phi) \cdot \exp \left(\phi(\gamma \ \lambda)  \right).
$$

 \ Thus,

$$
 {\bf E} e^{(\lambda,\xi)}  \le 1 + \prod_{k=1}^d |\lambda_k| \cdot
 I_{\gamma}(\phi) \cdot \exp \left(\phi(\lambda\gamma)  \right) \le
\exp \left(\phi(C \ \gamma \ \lambda)  \right),
$$
as long as $  \vec{\lambda} >> \vec{1}. $\par

 \  Another details are simple and may be omitted. \par

\vspace{4mm}

 \ As a slight consequence:\par

\vspace{4mm}

 \ {\bf  Corollary 2.1. } Suppose  as above that the function $  \phi(\cdot)  $ satisfies the conditions of theorem 2.1. We assert:
the centered non-zero random vector $  \xi  $ belongs to the space  $   B(\phi): $

$$
\exists C_1 \in (0,\infty), \ \forall \lambda \in R^d
 \Rightarrow  {\bf E} e^{(\lambda,\xi)} \le e^{ \phi( C_1 \cdot \lambda) }, \ \lambda \in R^d
$$
if and only if

$$
\exists C_2 \in (0,\infty), \ \forall \ x \in R^d_+  \Rightarrow \  U(\vec{\xi}, \vec{x}) \le \exp \left( - \phi^*(\vec{x}/C_2) \right).
$$

 \ More precisely, the following implication holds: there is certain finite positive constant $  C_3 =  C_3(\phi)  $ such that
  for arbitrary non-zero centered random vector $  \xi: \ ||\xi|| = ||\xi||B(\phi)  < \infty \ $ or equally

$$
 \forall \lambda \in R^d \Rightarrow  {\bf E} e^{(\lambda,\xi)} \le e^{ \phi( ||\xi|| \cdot \lambda) }
$$
iff

$$
\exists C_3(\phi) \in (0,\infty) \ \forall \ x \in R^d_+  \Rightarrow \  U(\vec{\xi}, \vec{x}) \le
\exp \left( - \phi^*(\vec{x}/(C_3 / ||\xi||) \right). \eqno(2.9)
$$

\vspace{3mm}

\ {\bf  Corollary 2.2. } Assume the non-zero centered random vector $  \xi = \{\xi(1), \xi(2), \ldots, \xi(d) \} $ belongs
to the space $  B(\phi):  $

$$
{\bf E} e^{(\lambda,\xi)} \le e^{ \phi(||\xi|| \cdot \lambda) }, \ \phi \in Y(R^d),
\eqno(2.10)
$$
and let $  y  $ be arbitrary positive non-random number. Then $  \forall y  > 0 \ \Rightarrow  $

$$
{\bf P} \left( \min_{j = 1,2,\ldots,n} |\xi(j)| > y  \right) \le 2^d \cdot
\exp \left(- \phi^*(y/||\xi||,y/||\xi||, \ldots, y/||\xi||) \right). \eqno(2.11)
$$

 \ The last estimate plays a very important role in the theory of discontinuous random fields, in particular, in the theory of
Central Limit Theorem in the space of Prokhorov-Skorokhod, see \cite{Grigorjeva1}, \cite{Ostrovsky5} - \cite{Ostrovsky7}.\par

\vspace{3mm}

{\bf Example 2.1.} Let as before $  \phi(\lambda) = \phi^{(B)}(\lambda)  = 0.5(B\lambda,\lambda),  $ where
$  B  $  is non-degenerate positive definite symmetrical matrix, in particular $  \det B > 0. $ It follows from theorem 2.1
that the (centered) random vector $ \xi $ is subgaussian relative the matrix $ B: $

$$
\forall \lambda \in R^d \ \Rightarrow {\bf E} e^{ (\lambda, \xi)} \le e^{0.5 (B \lambda, \lambda) ||\xi||^2 }.
$$
iff for some finite positive constant $  K = K(B,d) $ and for any  non-random positive vector $  x = \vec{x} $

$$
U(\xi,x) \le e^{- 0.5 \ \left( (B^{-1}x,x)/(K||\xi||^2) \right) }. \eqno(2.12)
$$

\vspace{3mm}

\ {\bf  Remark 2.1. } It is no hard to verify that the condition (2.5) follows immediately from the restriction (1.6). \par

\vspace{3mm}

\ {\bf  Remark 2.2. } The belonging of the random vector $ \xi = \vec{\xi} $ to the Banach  space $ B(\phi) $ may be characterized
as follows. Define the following {\it exponential} Young-Orlicz function

$$
N_{\phi}(u) = e^{ \phi^*(u) } - 1, \eqno(2.13)
$$
and introduce the Orlicz's space  $ L_{N} = L_{N_{\phi}}  $ over our source probability space
$ (\Omega,F,{\bf P} ) $  with correspondent $  N \ -  $ function $ N_{\phi}(u) . $  It is proved in \ \cite{Ostrovsky5} \ that

$$
\vec{\xi} \in L_{N_{\phi}}  \ \Longleftrightarrow
\exists C \in (0, \infty), \ U(\vec{\xi}, \vec{x}) \le \exp \left(- \phi^*(\vec{x}/C) \right). \eqno(2.14)
$$
 \ The one-dimensional case is provided in  \cite{Kozachenko1}. \par

\vspace{3mm}

\ {\bf  Remark 2.3. } The {\it exponential} exactness of the estimates (2.6) and (2.7) take place still in the
one - dimensional case, see  \ \cite{Kozachenko1}, \cite{Ostrovsky1}, chapter 1. \par

\vspace{3mm}

\ {\bf  Remark 2.4.} The conditions of theorem 2.1 are satisfied if for example the generating function $  \phi(\lambda) $ is
twice continuous differentiable regular varying function as infinity on the Euclidean norm $ |\lambda|, $ i.e. is radial, or
spherical function, with degree greatest than one:

$$
\phi(\lambda) = |\lambda|^{\beta} \ M(|\lambda|), \ |\lambda| \ge 1; \ \beta = \const > 1, \eqno(2.15)
$$
where $  M(z) $ is slowly varying as infinity  twice continuous differentiable function:

$$
\forall t > 0 \ \Rightarrow  \lim_{z \to \infty} \frac{M(t z)}{M(z)} = 1.
$$

 \ If in particular $ \ \phi(\lambda) = |\lambda|^{\beta}/\beta, \ |\lambda| \ge 1, \beta = \const > 1, $ then

$$
\phi^*(x) = \frac{\beta - 1}{\beta} \ |x|^{\beta/(\beta - 1)}, \ |x| \ge C(\beta). \eqno(2.16)
$$

\vspace{3mm}

 \ Of course, the classical Euclidean norm $  |\lambda| $  in (2.15) may be replaced on the other complete norm on the
space $  R^d. $ \par

\vspace{4mm}

\section{ Relation with moments. } \par

\vspace{4mm}

 \hspace{5mm} We intend in this section to simplify the known results and proofs for the moment estimates for the
one - dimensional r.v., see \cite{Kozachenko1},  \cite{Ostrovsky1}, p. 50-53,  and extend obtained result on the
multivariate case. \par
 \ We will use the following elementary inequality

$$
x^r \le \left( \frac{r}{\lambda \ e} \right)^r \cdot e^{\lambda \ x}, \ r,\lambda,x > 0, \eqno(3.0)
$$
and hence

$$
|x|^r \le \left( \frac{r}{\lambda \ e} \right)^r \ \cosh (\lambda x), \  \ r,\lambda > 0, \ x \in R. \eqno(3.0a)
$$

 \ As a consequence: let $  \xi $ be non-zero  one-dimensional mean zero random variable belonging to the space $  B(\phi).  $
Then

$$
{\bf E} |\xi|^r \le 2 \  \left( \frac{r}{\lambda \ e} \right)^r  \ e^{ \phi(\lambda ||\xi||) }, \ \lambda > 0. \eqno(3.1)
$$

  \  Authors of the article \cite{Kozachenko1}, see also \cite{Ostrovsky1}, chapter 1, section 1.5 have long chosen in the
inequality (3.1) the value

$$
\lambda = \lambda_0 := \phi^{-1}(r/||\xi||).
$$

 \ We intend here to choose the value $ \lambda $ for reasons of optimality.
 \ We hope that this method is more simple and allows easy multivariate generalization. \par

 \ In detail, introduce the function

$$
\Phi(\mu) := \phi^* \left(e^{\mu} \right), \ \mu \in R. \eqno(3.2)
$$

 \ Suppose temporarily  for simplicity in (3.1)  $ \ ||\xi|| = ||\xi||B(\phi) = 1. $ One can rewrite  (3.1) as follows.

$$
{\bf E} |\xi|^r \le 2 \ r^r \ e^{-r} \ e^{ - r \ln \lambda + \phi(\lambda ) } =
$$

$$
 2 \ r^r \ e^{-r} \ e^{ - r \mu + \phi \left(e^{\mu} \right) } = 2 \ r^r \ e^{-r} \ e^{ -\left(r \mu - \Phi(\mu) \right) }, \eqno(3.3)
$$
and we deduce after minimization over $ \mu \ $ (or equally over $ \lambda)  $

$$
{\bf E} |\xi|^r \le  2 \ r^r \ e^{-r} \ e^{- \Phi^*(r) }.
$$

 \ So, we proved in fact the following statement. \par

 \vspace{3mm}

{\bf Proposition 3.2.} Let $ \phi (\cdot) $  be arbitrary non-negative continuous function and let the centered numerical
r.v. $  \xi $  be such that $  \xi \in B(\phi) $ or equally

$$
U(\xi,x) \le \exp(-\phi^*(x)), x \ge 0.
$$

Then

$$
|\xi|_r \le 2^{1/r} \ r \ e^{-1} \  e^{- \Phi^*(r)/r } \ ||\xi||B(\phi), \ r > 0. \eqno(3.4)
$$

 \ The inverse conclusion as well as the multivariate generalization contains in the preprint \cite{Ostrovsky5}. {\it But
we will demonstrate further a more effective approach.} \par

\vspace{3mm}

\ Let again \ (temporarily) $ \ d = 1 $ and suppose $  \xi \in B(\phi), \ ||\xi|| := ||\xi||B(\phi) \in (0,\infty). $ Then

$$
{\bf E} |\xi|^p \le 2 p \int_0^{\infty} x^{p-1} e^{ - \phi^*(x/||\xi||) } \ dx  =
2 \ ||\xi||^p \int_{-\infty}^{\infty} \exp \left( pz - \phi^*(e^z) \right) \ dz=
$$

$$
2 \ ||\xi||^p \ \int_{-\infty}^{\infty} \exp \left( pz - \Phi(z) \right) \ dz. \eqno(3.5)
$$

 \ We deduce using the saddle-point method:

$$
{\bf E} |\xi|^p \le 2 \ ||\xi||^p  \cdot \exp \left( \sup_p ( (C(\phi) \cdot p) \ z - \Phi(z)   )    \right) =
$$

$$
2 \ ||\xi||^p \ \exp \left(  \Phi^*(C(\phi) \ p)  \right), \ p \ge 1,
$$
or equally

$$
|\xi|_p \le
 2^{1/p} \ ||\xi|| \ \exp \left(  \Phi^*(C(\phi) \ p)/p  \right) \le
$$

$$
 ||\xi|| \ \exp \left(  \Phi^*(C_1(\phi) \ p)/p  \right), \ p \ge 1. \eqno(3.6)
$$

 \ The rigorous proof will be carried out later, in the multidimensional case. \par

 \ Conversely, let the inequality (3.6) be given. We can suppose for simplicity

$$
|\xi|_p  \le  \exp \left(  \Phi^*(p)/p  \right), \ p \ge 1. \eqno(3.7)
$$

 \ Then

$$
{\bf E} |\xi|^p \le \exp \left(  \Phi^*(p) \right).
$$
 \ Tchebychev's inequality gives us

$$
{\bf P}(|\xi| > u) \le \frac{ \exp \left(  \Phi^*(p) \right)}{ u^p} =
\exp \left( -(p \ \ln u - \Phi^*(p))  \right), \ u \ge 2,
$$
therefore

$$
{\bf P}(|\xi| > u) \le \exp \left( - \sup_p(p \ \ln u - \Phi^*(p))  \right) =  \exp \left( - \Phi^{**}(\ln u) \right) =
$$

$$
\exp \left( - \Phi(\ln u)   \right) = \exp \left( - \phi^*(u)  \right)
$$
by virtue of theorem Fenchel-Moraux. \par

\vspace{3mm}

 {\bf Example 3.1.}  Suppose that the function $ \phi (\cdot)  $  is  such that for some
 constant $  p  > 1 $

$$
\phi(\lambda) = \phi_p(\lambda) \le C_1 \ |\lambda|^p, \ |\lambda| >1.
$$
 \ Let  also the centered non-zero random variable  $  \xi $ belongs to the space $ G\psi_p.  $
Then

$$
|\xi|_r \le C_2(p) \ r^{1/q} \ ||\xi||B(\phi_r), \ r \ge 1, \ q := p/(p-1),
$$
and the inverse conclusion is also true: if $ {\bf E} \xi = 0 $ and if for some positive finite constant $  K $

$$
\forall r \ge 1 \ \Rightarrow  |\xi|_r \le K \ r^{1/q},
$$
then $ \xi(\cdot) \in B(\phi_p) $  and wherein $ ||\xi||  B(\phi_p) \le C_3 \ K.   $ \par

\vspace{4mm}

 {\bf Example 3.2.}  Suppose now that the function $ \phi (\cdot)  $ is  such that

$$
\phi(\lambda) = \phi^{(K)}(\lambda)  \le \frac{ C_4}{K - \ |\lambda|}, \ |\lambda| < K.
$$
 \ Let  also the centered non-zero random variable  $  \xi $ belongs to this space $ G\psi^{(K)}.  $
Then

$$
|\xi|_r \le C_5 \ K \ r \ ||\xi||B(\phi^{(K)}),
$$
and likewise the inverse conclusion is also true: if $ {\bf E} \xi = 0 $ and if for some finite positive constant $  K $

$$
\forall r \ge 1 \ \Rightarrow  |\xi|_r \le K \ r,
$$
then $ \xi(\cdot) \in B(\phi^{(K)}) $  and herewith $ ||\xi|| B(\phi^{(K)}) \le C_6 \ K.   $ \par

\vspace{4mm}

 \ We need  getting to the presentation of the multidimensional case to extend our notations and restrictions.
 In what follows in this section the variables $ \lambda, r, x, \xi  $  are as before vectors from the space $  R^d, \ d = 2,3,\ldots, $
and besides $  r = \vec{r} = \{ r(1), r(2), \ldots, r(d) \}, \ r(j) \ge 1.  $\par

 \ Vector notations:

 $$
  |r| = |\vec{r}| = \sum_j r(j), \hspace{4mm} |\xi| = |\vec{\xi}| = \{ |\xi(1)|, |\xi(2)|, \ldots, |\xi(d)| \} \in R^d_+,
 $$

$$
 \vec{x} \ge 0 \ \Leftrightarrow \forall j \hspace{3mm} x(j) \ge 0;
$$

 $$
 x^r = \vec{x}^{\vec{r}} = \prod_{j=1}^d x(j)^{r(j)}, \ \vec{x} \ge 0,
 $$

$$
\ln \vec{\lambda} = \{  \ln \lambda(1), \ \ln \lambda(2), \ldots, \ \ln \lambda(d) \}, \hspace{3mm} \vec{\lambda} > 0,
$$

$$
e^{\vec{\mu}} = \{ e^{\mu(1)}, \ e^{\mu(2)}, \ldots, \  e^{\mu(d)} \},
$$

$$
\Phi(\mu) = \Phi(\vec{\mu}) = \phi^* \left(e^{\vec{\mu}} \right),
$$

 $$
 \frac{r}{\lambda \cdot  e}  =  \frac{\vec{r}}{ \vec{\lambda} \ \cdot e} = \prod_{j=1}^d \left( \frac{r(j)}{ e \ \lambda(j)} \right) =
 e^{-|r|} \cdot \prod_{j=1}^d \left( \frac{r(j)}{\lambda(j)} \right),
 $$

$$
|\xi|_r = |\vec{\xi}|_{\vec{r}} = \left( {\bf E}|\vec{\xi}|^{\vec{r}} \right)^{1/|r|}.
$$

 \ We will use again the following elementary inequality

$$
x^r \le \left( \frac{r}{\lambda \ e} \right)^r \cdot e^{(\lambda, \ x)}, \hspace{3mm} r, \ \lambda, \ x > 0. \eqno(3.13)
$$

 \ As a consequence: let $  \xi $ be non-zero  $ d  - $  dimensional mean zero random vector belonging to the space $  B(\phi).  $
Then

$$
{\bf E} |\xi|^r \le 2^d \  \left( \frac{r}{\lambda \ e} \right)^r  \ e^{ \phi(\lambda ||\xi||) } =
2^d \ e^{-|r|} \ r^r \ \lambda^{-r} \ e^{ \phi(\lambda ||\xi||) }, \ \lambda > 0. \eqno(3.14)
$$

 \ We find likewise to the one - dimensional case: \\

{\bf Proposition 3.3.} Let $ \phi (\cdot) $  be arbitrary non-negative continuous function and let the centered numerical
r.v. $  \xi $  be such that $  \xi \in B(\phi): \ 0 < ||\xi|| = ||\xi||B(\phi) < \infty. $ \par

 \ Then

$$
| \vec{\xi}|_{\vec{r}} \le e^{-1} \cdot  2^{d/|r|} \cdot \prod_j r(j)^{r(j)/|r|}  \cdot  e^{- \Phi^*(r)/|r| } \cdot ||\xi||B(\phi),
 \ r = \vec{r} > 0. \eqno(3.15)
$$

 \ Note that in general case the expression $ |\xi|_r  $ does not represent the norm relative the random vector $  \vec{\xi}. $
 But if we denote

$$
\psi_{\Phi}(\vec{r}) :=  e^{-1} \cdot 2^{d/|r|} \cdot \prod_j r(j)^{r(j)/|r|}  \cdot  e^{- \Phi^*(r)/|r| }
$$
  and define
$$
||\xi||G\psi_{\Phi} \stackrel{def}{=} \sup_{\vec{r} \ge 1} \left[  \frac{| \vec{\xi}|_{\vec{r}} }{\psi_{\Phi}(\vec{r}) }  \right],
$$
we obtain some generalization of the known one-dimensional Grand Lebesgue Space (GLS) norm, see \cite{Kozachenko1},
\cite{Fiorenza1}-\cite{Fiorenza3}, \cite{Iwaniec1}-\cite{Iwaniec2}, \cite{Ostrovsky1}, chapter 1. \par
 \ These multivariate  generalization of a form $ ||\xi||G\psi_{\Phi} $ based in turn on the theory of the so-called
mixed (anisotropic) Lebesgue-Riesz spaces \cite{Besov1}, chapters 1,2; appears at first perhaps in the authors preprints
\cite{Ostrovsky8}-\cite{Ostrovsky9}. \par

\vspace{3mm}

 \ The statement of proposition 3.3 may be rewritten as follows.

$$
||\xi||G\psi_{\Phi}  \le ||\xi||B(\phi). \eqno(3.15a)
$$

 \ Let us state the inverse up to multiplicative constant inequality. \par

\vspace{4mm}

{\bf Theorem 3.1.} Suppose $  \phi(\cdot) $  is  Young-Orlicz function.  Let $ {\bf E} \xi =0  $ and let
$ K = ||\xi||G\psi_{\Phi} < \infty. $  Suppose also that there is a constant $  \gamma \in (0,1)  $ such that
the following integral converges:

$$
L := \int_{R^d} e^{\Phi(\gamma z) - \Phi(z) } \ dz < \infty. \eqno(3.16)
$$

 \ Then $ \xi \in B(\phi)  $ and moreover both the norms
$ || \xi||B(\phi) $  and  $ ||\xi||G\psi_{\Phi} $ are equivalent:

$$
 ||\xi||G\psi_{\Phi} \le || \xi||B(\phi) \le C_4(\phi) \ ||\xi||G\psi_{\Phi}, \   C_{4}(\phi) \in (0,\infty).  \eqno(3.17)
$$

\vspace{3mm}

 \ The {\bf proof} of this theorem is at the same as one in the theorem  2.1.  Let at first $ 0 < ||\xi||G\psi_{\Phi} < \infty;  $ we
can  agree without loss of generality

$$
| \ \vec{\xi} \ |_{\vec{p}}^{ |p| } = {\bf E} |\xi|^{|p|}
\le \exp \left( \Phi^*(\vec{p})  \right).
$$
 \ We apply again the Tchebychev-Markov inequality:

$$
U(\vec{\xi}, \vec{x}) \le \exp \left \{- \left( \sum_j p_j \ \ln x(j) - \Phi^*(\vec{p}) \right) \right\}, \ \vec{x} \ge 2 \ \vec{1},
$$
 or after the minimization over $  \vec{p}  $

$$
U(\vec{\xi}, \vec{x}) \le  \exp \left\{  -\sup_{ \vec{p}} \left( \sum_j p(j) \ \ln x(j) - \Phi^*(\vec{p} ) \right) \right\} =
$$

$$
\exp \left\{ -\Phi^{**}(\ln \vec{x}) \right\} =
\exp \left\{ -\Phi(\ln \vec{x}) \right\} = \exp \left\{ -\phi^*(\vec{x}) \right\}
$$
again by virtue of theorem Fenchel-Moraux. \par
 \ The case $ \exists k \ \Rightarrow x(k) < 2 $ may be easily reduced  to the considered here. \par

\vspace{3mm}

{\bf Remark 3.1.} Note that we do not use in this direction the condition (3.16). \par

\vspace{3mm}

 \ Let us prove the inverse conclusion. Suppose the mean zero random vector $ \xi = \vec{\xi}  $ is such that
 $ ||\xi|| = ||\xi||B(\phi) = 1.  $ We intend to estimate the (power) moment $  |\xi|_p^{|p|} = {\bf E}|\vec{\xi}|^{ \vec{p} } =
 {\bf E}|\xi|^{p}. $ We deduce using the tail estimate (Theorem 2.1) and after integration by parts

$$
 {\bf E}|\xi|^{p} \le 2^d \cdot \prod_j p_j \cdot \int_{R^d_+} \prod_j x(j)^{ p(j) - 1} \ e^{ -\phi^*(x)} \ dx =
$$

$$
2^d \ \int_{R^d} e^{ (p,z) - \phi^*(\exp z) } \ dz = 2^d \ \int_{R^d} e^{ (p,z) - \Phi(z)  } \ dz \stackrel{def}{=} 2^d \ Q(p).\eqno(3.18)
$$

 \ Let $  \gamma = \const \in (0,1).  $   We can write acting as before, in the second section,

$$
(p,z) \le \Phi(\gamma z) + \Phi^*(p /\gamma).
$$

  We have after substituting into (3.18)

$$
Q(p) \le \exp \left(\Phi^*(p/\gamma)  \right) \cdot \int_{R^d} \exp \left( \Phi(\gamma x) - \Phi(x)  \right) \ dx =
$$

$$
L \cdot\exp \left(\Phi^*(p/\gamma)  \right).\eqno(3.19)
$$
 \ This completes the proof of theorem 3.1. \par

 \vspace{3mm}

{\sc  To summarize.}  For arbitrary centered $ d \ -  $ dimensional random vector $ \vec{\xi}  $ under formulated above, in
theorems 2.1 and 3.1 conditions the following predicated are equivalent:

$$
{\bf A. } \hspace{5mm} \vec{\xi} \in L_{N_{\phi}}. \eqno(3.20a)
$$

$$
{\bf B. }  \hspace{5mm} \exists C \in (0, \infty), \ U(\vec{\xi}, \vec{x}) \le \exp \left(- \phi^*(\vec{x}/C) \right). \eqno(3.20b)
$$

$$
{\bf C.}  \hspace{5mm} \vec{\xi} \in G\psi_{\Phi}. \eqno(3.20c)
$$

$$
{\bf D.}  \hspace{5mm} \vec{\xi} \in B(\phi). \eqno(3.20d)
$$

 \ Herewith  all the  Banach norms $ ||\xi||B(\phi), \  ||\xi|| L_{N_{\phi}} $ and $ ||\xi|| G\psi_{\Phi}  $ are
linear equivalent. \par

 \vspace{3mm}

\section{Concluding remarks.}

\vspace{3mm}

{\bf A.} To mention the possible application of obtained results. At first: a non-asymptotical {\it exponential} bounds
for the normed sums of random vectors, in the spirit of the classical article of Yu.V.Prokhorov \  \cite{Prokhorov1}, see
\ \cite{Ostrovsky5}.  Also one can derive an asymptotical and non-asymptotical analysis of  discontinuous random fields, to
generalize the main result of \ \cite{Ostrovsky6}. \par
 \ Another applications-an investigation of the Central Limit Theorem for discontinuous random fields, \ \cite{Ostrovsky7},
with further applications in the non-parametrical statistics as well as in the Monte-Carlo method \ \cite{Grigorjeva1}. \par

\vspace{3mm}

{\bf B.} It may be studied easily the case when the our considered generating Young-Orlicz function
$  \phi = \phi(\lambda), \ \lambda \in R^d $ is finite only in some convex non-trivial neighborhood $ V $ of origin,
likewise the one-dimensional case, see e.g. \cite{Kozachenko1},  \cite{Ostrovsky1}, chapter 1. \par

\vspace{3mm}

{\bf C.} The another (asymptotical) approach  for research of the non-linear functionals $  F(\xi) = F(\xi(\cdot)) $ of the random fields
may be  found in the famous survey of V.I.Piterbarg and V.R.Fatalov \ \cite{Piterbarg1}. This approach based on the investigation of
the  asymptotical behavior as $  \lambda \to 0+ $ for the {\it Laplace transform} from  the considered  functional

$$
\hat{L}_F(\lambda) := {\bf E} \ e^{ - \lambda F(\xi) }.
$$

\end{document}